\documentclass{article}

\usepackage{amsmath}
\usepackage{amssymb}
\usepackage{graphicx}
\usepackage{bm}
\usepackage{mathdots}
\usepackage{booktabs}
\usepackage{rotating}
\usepackage{booktabs, multirow, array}
\usepackage{float}
\usepackage[ruled,linesnumbered]{algorithm2e}
\usepackage[a4paper,left=2.8cm,right=2.8cm,top=2.5cm,bottom=2.5cm]{geometry}
\usepackage{fancyhdr}
\usepackage{hyperref}
\usepackage{multirow} 
\usepackage[framemethod=tikz]{mdframed}
\newtheorem{theorem}{Theorem}[section]
\newtheorem{corollary}{Corollary}[section]

\newtheorem{lemma}{Lemma}[section]
\newtheorem{definition}{Definition}[section]
\newtheorem{assumption}{Assumption}[section]
\newtheorem{remark}{Remark}[section]

\usepackage[mathscr]{euscript}

\newcommand{\measure}{\text{d}\omega}

\makeatletter % `@' now normal "letter"
\@addtoreset{equation}{section}
\makeatother  % `@' is restored as "non-letter"

\newenvironment{proof}{{\noindent\it Proof.}\quad}{\hfill $\square$\\}
\usepackage{url}
\usepackage{mathtools}

\begin{document}
\title{Recovery Thresholding Hyperinterpolations in Signal Processing}

\author{Congpei An\footnotemark[1]
       \quad \text{and}\quad Jiashu Ran\footnotemark[2]}
\renewcommand{\thefootnote}{\fnsymbol{footnote}}
\footnotetext[1]{School of Mathematics and Statistics, Guizhou University, Guiyang 550025, Guizhou, China (ancp@gzu.edu.cn,andbachcp@gmail.com),}
\footnotetext[2]{Department of Mathematical and Statistical Sciences, University of Alberta, Edmonton, Alberta T6G 2G1, Canada (jran@ualberta.ca)}
\maketitle

\begin{abstract}
This paper introduces recovery thresholding hyperinterpolations, a novel class of methods for sparse signal reconstruction in the presence of noise. We develop a framework that integrates thresholding operators—including hard thresholding, springback, and Newton thresholding—directly into the hyperinterpolation structure to maintain sparsity during signal recovery. Our approach leverages Newton's method to minimize one-dimensional nonconvex functions, which we then extend to solve multivariable nonconvex regularization problems. The proposed methods demonstrate robust performance in reconstructing signals corrupted by both Gaussian and impulse noise. Through numerical experiments, we validate the effectiveness of these recovery thresholding hyperinterpolations for signal reconstruction and function denoising applications, showing their advantages over traditional approaches in preserving signal sparsity while achieving accurate recovery.
\end{abstract}



\section{Introduction}
Modern signal processing and data analysis face unprecedented challenges due to the increasing complexity and volume of data across diverse applications. While convex optimization methods have been the cornerstone of many data recovery techniques, they often prove inadequate when dealing with highly corrupted signals or when seeking sparse solutions. This limitation has motivated the development of nonconvex optimization approaches \cite{an2022springback,Chen2014complexity_nonconvex,Chen2010lower_bound}, which offer enhanced recovery capabilities at the cost of increased computational complexity.

In this paper, we introduce recovery thresholding hyperinterpolations—a novel class of methods that bridges the gap between robust signal recovery and computational efficiency. Our approach builds upon the hyperinterpolation framework established by Sloan \cite{sloan1995hyperinterpolation}, which uses polynomial approximations to reconstruct functions from discrete observations. By incorporating thresholding operators \cite{donoho1994ideal}  directly into the hyperinterpolation structure, we create methods that naturally enforce sparsity while maintaining the approximation power.

The key innovation of our work lies in the systematic integration of various thresholding strategies within the hyperinterpolation framework. We develop and analyze three distinct variants: hard thresholding hyperinterpolation \cite{an2023hard}, springback hyperinterpolation \cite{an2022springback}, and Newton hyperinterpolation. Each variant corresponds to a specific nonconvex regularization scheme, offering different trade-offs between recovery accuracy and computational complexity. The part of sparsity analysis uses concentration inequalities, providing probabilistic guarantees on the support recovery.

The practical significance of our methods is demonstrated through applications to signal reconstruction and function denoising in the presence of both Gaussian and impulse noise. Our numerical experiments show that recovery thresholding hyperinterpolations outperform traditional approaches in preserving signal structure while effectively suppressing noise.

The remainder of this paper is organized as follows. Section 2 provides the necessary background on hyperinterpolation and quadrature rules. Section 3 presents our analysis of one-dimensional nonconvex functions and derives the optimal thresholding operators. Section 4 extends these results to multivariable nonconvex regularization problems. Section 5 develops the sparsity analysis framework using concentration inequalities. Section 6 presents numerical experiments validating our theoretical findings.

\section{Backgrounds on hyperinterpolation}\label{sec:setting}
Let $\Omega$ be a compact set of $\mathbb{R}^s$ with a positive measure $\omega$.  Suppose $\Omega$ has finite measure with respect to $\measure$, that is,
\begin{eqnarray*}
    \int_{\Omega} \,\rm{d} \omega = V < \infty.
\end{eqnarray*}

We denote by $L_2(\Omega)$ the Hilbert space of square-integrable functions on $\Omega$ with the $L_2$ inner product
\begin{eqnarray}\label{equ:continuousinner}
\langle f ,g \rangle = \int_{\Omega} fg \,{\rm{d}} \omega  \qquad \forall f,g \in L_2(\Omega),
\end{eqnarray}
and the induced norm $\|f\|_2: = \langle f,f\rangle^{1/2}$. Let $\Pi_n (\Omega) \subset L_2(\Omega)$ be the linear space of polynomials with total degree not strictly greater than $n$, restricted to $\Omega$, and let $d_n= \dim (\Pi_n(\Omega))$ be the dimension of $\Pi_n(\Omega)$.

Next we define an orthonormal basis of $\Pi_n(\Omega)$
%\begin{eqnarray}\label{set:orthonormalbasis}
\[\{ \Phi_{\ell}| \ell = 1, \ldots, d_n \} \subset \Pi_n(\Omega)\]
%\end{eqnarray}
in the sense of
%\begin{eqnarray}\label{equ:kronecker}
\[\langle \Phi_{\ell},\Phi_{\ell'} \rangle = \delta_{\ell \ell'} \qquad \forall 1 \leq \ell, \ell' \leq d_n.\]
%\end{eqnarray}

The $L_2(\Omega)$-\emph{orthogonal projection} $\mathcal{T}_n: L_2(\Omega) \to \Pi_n(\Omega)$ can be uniquely defined by
\begin{equation}\label{def:L2Projection}
\mathcal{T}_n  f:= \sum_{\ell=1}^{d_n} \hat{f}_{\ell} \Phi_{\ell} = \sum_{\ell=1}^{d_n} \langle f, \Phi_{\ell} \rangle \Phi_{\ell}  , \qquad \forall f \in L_2(\Omega),
\end{equation}
where $\{\hat{f}_{\ell}\}_{\ell=1}^{d_n}$ are the Fourier coefficients
\begin{eqnarray*}
    \hat{f}_{\ell} := \langle f, \Phi_{\ell} \rangle = \int_{\Omega} f\Phi_{\ell} \,{\rm{d}} \omega, \qquad \forall \ell = 1,\ldots,d_n.
\end{eqnarray*}

In order to numerically evaluate the scalar product in \eqref{def:L2Projection}, it is fundamental to consider an $N$-point quadrature rule of PI-type (Positive weights and Interior nodes), i.e.,
\begin{eqnarray}\label{equ:Npointquad}
\sum_{j=1}^{N} w_{j}g(\mathbf{x}_j) \approx \int_{\Omega} g \,{\rm{d}} \omega \qquad \forall g \in \mathcal{C}(\Omega),
\end{eqnarray}
where the quadrature points $\{\mathbf{x}_1, \ldots, \mathbf{x}_N\}$ belong to $\Omega$ and the corresponding quadrature weights $\{w_1, \ldots, w_N\}$ are positive, and $\mathcal{C}(\Omega)$ is a continuous function space. Furthermore we say that \eqref{equ:Npointquad} has algebraic degree of exactness $\delta$ if
\begin{eqnarray}\label{equ:exactness}
\sum_{j=1}^{N} w_j p(\mathbf{x}_j) = \int_{\Omega}p \,{\rm{d}} \omega    \qquad \forall p \in\Pi_{\delta}(\Omega).
\end{eqnarray}
With the help of a quadrature rule \eqref{equ:Npointquad} with algebraic degree of exactness $\delta=2n$, we can introduce a ``\emph{discrete (semi) inner product}'' on  $\mathcal{C}(\Omega)$ \cite{sloan1995hyperinterpolation} by
\begin{eqnarray}\label{equ:discreteinner}
\langle f,g \rangle_N := \sum_{j=1}^{N} w_j f(\mathbf{x}_j) g(\mathbf{x}_j) \qquad \forall f,g \in \mathcal{C}(\Omega),
\end{eqnarray}
corresponding to the $L_2(\Omega)$-inner product \eqref{equ:continuousinner}. For any $p,q \in \Pi_{n}(\Omega)$, the product $pq$ is a polynomial in $\Pi_{2n}(\Omega)$. Therefore it follows from the quadrature exactness of \eqref{equ:exactness} for polynomials of degree at most $2n$ that
%\begin{eqnarray}\label{equ:kroneckerdiscrete}
\[\langle {p,q} \rangle_N = \langle {p,q} \rangle = \int_{\Omega} pq \,{\rm{d}} \omega, \qquad \forall p, q \in \Pi_{n}(\Omega).\]
%\end{eqnarray}

In 1995, Sloan introduced in \cite{sloan1995hyperinterpolation} the \emph{hyperinterpolation operator} $\mathcal{L}_n: \mathcal{C}(\Omega) \to \Pi_n(\Omega)$ as
\begin{eqnarray}\label{equ:hyper}
\mathcal{L}_n f:= \sum_{\ell=1}^{d_n} \langle {f, \Phi_{\ell}} \rangle_N   \Phi_{\ell}.
\end{eqnarray}
$\mathcal{L}_n f$ is a projection of $f$ onto $\Pi_n(\Omega)$ obtained by replacing the $L_2(\Omega)$-inner products \eqref{equ:continuousinner} in the $L_2(\Omega)$-orthogonal projection $\mathcal{T}_n f$ by the discrete inner products \eqref{equ:discreteinner}.

To explore one of its important features, we introduce the following discrete least squares approximation problem
\begin{eqnarray}\label{HyperLeastApprox}
\min_{p \in \Pi_n(\Omega) }  \left\{ \frac{1}{2} \sum_{j=1}^{N}w_j[p(\mathbf{x}_j) - f(\mathbf{x}_j)]^2 \right\}
\end{eqnarray}
with $p(\mathbf{x}) = \sum_{\ell=1}^{d_n}\alpha_{\ell}\Phi_{\ell}(\mathbf{x}) \in \Pi_{n}(\Omega)$, or equivalently
%\begin{eqnarray}\label{HyperLeastApproxMatrix}
\[\min\limits_{\bm{\alpha} \in \mathbb{R}^{d_n}}~~ \frac{1}{2} \|\mathbf{W}^{1/2} ( \mathbf{A}\bm{\alpha} - \mathbf{f})\|_2^2,\]
%\end{eqnarray}
where $\mathbf{W}=\text{diag}(w_1,\ldots,w_N)$ is the quadrature weights matrix, $\mathbf{A}=(\Phi_{\ell}(\mathbf{x}_j))\in\mathbb{R}^{N\times {d_n}}$ is the sampling matrix, $\bm{\alpha}=[\alpha_1,\ldots,\alpha_{d_n}]^{\text{T}}\in\mathbb{R}^{d_n}$ and $\mathbf{f}=[f(\mathbf{x}_1),\ldots,f(\mathbf{x}_N)]^{\text{T}}\in\mathbb{R}^N$ are two column vectors (recall $\mathbf{x}_j\in\mathbb{R}^s$). Sloan in \cite{sloan1995hyperinterpolation}  revealed that the relation between the hyperinterpolation $\mathcal{L}_n f$ and the best discrete least squares approximation (weighted by quadrature weights) of $f$ at the quadrature points. More precisely he proved the following important result:
\begin{lemma}[Lemma 5 in \cite{sloan1995hyperinterpolation}]\label{prop:HyperLeastApprox}
Given $f\in\mathcal{C}(\Omega)$, let $\mathcal{L}_n{f}\in\Pi_n(\Omega)$ be defined by \eqref{equ:hyper}, where the discrete scalar product $\langle f,g \rangle_N$ in (\ref{equ:discreteinner}) is defined by an $N$-point quadrature rule of PI-type in $\Omega$ with algebraic degree of exactness $2n$. Then $\mathcal{L}_n{f}$ is the unique solution to the approximation problem \eqref{HyperLeastApprox}.
\end{lemma}

\begin{assumption}\label{assumption1}
The $N$-point quadrature rule \eqref{equ:Npointquad} of PI-type, with nodes $\mathbf{x}_j\in \Omega$ and weights $w_j>0$ for $j=1,2,\ldots,m$, has exactness degree $n+k$ with $0<k \leq n$, where $n,k \in \mathbb{N}$.
\end{assumption}

\begin{definition}\label{def:hyper_relax}\cite[Hyperinterpolation with an exactness-relaxing quadrature rule]{an2022quadrature}
Let $\langle \cdot, \cdot \rangle_N$ be an quadrature rule fulfilling Assumption \ref{assumption1}. Given $f \in \mathcal{C}(\Omega)$, the relaxed hyperinterpolant of degree $n$ to $f$ is defined as 
\begin{equation}
    \mathcal{R}_nf:=\sum_{\ell=1}^{d_n} \langle f, \Phi_{\ell} \rangle_N \Phi_{\ell}.
\end{equation}
\end{definition}

\begin{assumption}\label{assumption2}
We assume that there exists an $\eta \in [0,1)$, which is independent of $n$ and $p$, such that 
\begin{equation}
    \left| \sum_{j=1}^{N} w_jp(\mathbf{x}_j)^2 - \int_{\Omega} p^2 {\rm{d}}\omega \right| \leq \eta \int_{\Omega} p^2 {\rm{d}}\omega \quad \forall p \in \Pi_n{\Omega}.
\end{equation}
\end{assumption}

\begin{definition}\label{def:hyper_bypass}\cite[Hyperinterpolation bypass the quadrature exactness]{an2024bypassing}
Let $\langle \cdot, \cdot \rangle_N$ be an quadrature rule only satisfying  Assumption \ref{assumption2}. Given $f \in \mathcal{C}(\Omega)$, the unfettered hyperinterpolant of degree $n$ to $f$ is defined as 
\begin{equation}
    \mathcal{U}_n f:= \sum_{\ell=1}^{d_n} \langle f, \Phi_{\ell} \rangle_N \Phi_{\ell}.
\end{equation}
\end{definition}

\section{A Basic Analysis of the Minimizer for a One-Variable Non-convex Function}
Does $x$ have a closed form when $h(x) $ reaches its minimum?
\begin{equation*}
    h(x)= (y-x)^2+\lambda|x|^q,\quad q\in(0,1)\quad \lambda>0.
\end{equation*}

\begin{lemma}\label{lem:newton}
Let $ h(x)= (y-x)^2+\lambda|x|^q$ with $q\in(0,1)$ and $\lambda>0$. Then $h(x)$ reaches its minimum at 
%\begin{equation}\label{equ:newton}
%\psi_H(x;\lambda):=\left\{
%\begin{array}{ll}
%0, & \text{if } {\lambda^{\ast}}\leq \lambda, \\
%x', & \text{else},
%\end{array}
%\right.
%\end{equation}
\begin{equation}\label{equ:newton}
\psi_H(y;\lambda):=\left\{
\begin{array}{ll}
x', & \text{if } |y| > a , \\
0 , & \text{if } |y| \leq  a,
\end{array}
\right.
\end{equation}
where $x'$ is obtained via Newton's method by solving $h'(x)=0$ on $(\min\{0,y\},\max\{0,y\})$ and 
\[
a = \frac{2-q}{2} \cdot (1-q)^{\frac{q-1}{2-q}} \cdot \lambda^{\frac{1}{2-q}}
\]
%\[{\lambda^{\ast}} = \frac{4}{2^q(1-q)}\left( \frac{1-q}{2-q} |y|\right)^{2-q}.\]
In particular, we have 
\[
|x'| > \left(\frac{\lambda q (1-q)}{2}\right)^{\frac{1}{2-q}}.
\]

\end{lemma}
\begin{proof}
Without loss of generality, we assume $y>0$ since symmetric analysis applies, leading to analogous conclusions for $y<0$. Then 
\[h'(x)=\left\{\begin{array}{ll}
-2(y-x)-\lambda q |x|^{q-1}<0, & \text{if } x<0, \\
-2(y-x)+\lambda q x^{q-1}>0, & \text{if } x>y,
\end{array}
\right.
\]
which implies that
\[
\left\{\begin{array}{ll}
h(x)>h(0)=y^2, & \text{if } x<0, \\
h(x)>h(y)=\lambda y^q, & \text{if } x>y.
\end{array}
\right.
\]
For $x \in [0,y]$, we have 
\[
h'(x)=-2(y-x)+\lambda q x^{q-1},
\]
which implies that $h'(0^{+})=+\infty$ and $h'(y)=\lambda q y^{q-1}>0$. We claim that $h(x)$ can have at most two stationary points. 

Consider solving 
\[
h'(x)=-2(y-x)+\lambda q x^{q-1}=0, \qquad \forall x \in (0,y),
\]
which means that
\begin{equation}\label{equ:3-2}
    x^{1-q}(y-x)=\frac{\lambda q}{2} \qquad \forall x \in (0,y).
\end{equation}
Let $f(x)=x^{1-q}(y-x)$. We have 
\[
f'(x)=f(x)\left(\frac{1-q}{x} - \frac{1}{y-x} \right) \qquad \forall x \in (0,y).
\]
Noting that \[\frac{1-q}{x} - \frac{1}{y-x}\] is decreasing in $x$ and must have a root $\bar{x}$ on $(0,y)$. Indeed, since
\[
\frac{1-q}{\bar{x}} - \frac{1}{y-\bar{x}}=0 ,\]
we have
\[\bar{x} = \frac{1-q}{2-q}y \in (0,y).
\]
Then we obtain
\[
f(\bar{x})= \frac{1}{1-q}\left(\frac{1-q}{2-q}y\right)^{2-q}.
\]
If there exists $\bar{\lambda}$ such that $f(\bar{x})=\frac{\bar{\lambda} q}{2}$ in \eqref{equ:3-2}, then
\[
\bar{\lambda} = \frac{2}{q(1-q)}\left(\frac{1-q}{2-q}y\right)^{2-q}.
\]
Hence, the equation \eqref{equ:3-2} has two solutions if $\lambda < \bar{\lambda}$; $f(x)=\frac{\lambda q}{2}$ has one solution if $\lambda=\bar{\lambda}$; $f(x)=\frac{\lambda q}{2}$ has no solutions if $\lambda >\bar{\lambda}$. Therefore, $h(x)$ has two stationary points if $\lambda < \bar{\lambda}$; $h(x)$ has one stationary point if $\lambda=\bar{\lambda}$; $h(x)$ has no stationary points if $\lambda >\bar{\lambda}$.

Therefore, $h(x)$ can have at most two stationary points in $(0,y)$. Note that 
\[
\left\{\begin{array}{ll}
h(0)>h(y),  & \text{if } \lambda<y^{2-q}, \\
h(0)=h(y),  & \text{if } \lambda=y^{2-q}, \\
h(0)<h(y),  & \text{if } \lambda>y^{2-q}.
\end{array}\right.
\]
Here, we briefly analyze the relationship between $y^{2-q}$ and $\bar{\lambda}$. Note that $\frac{2}{q(1-q)}\left(\frac{1-q}{2-q}\right)^{2-q}$ decreases from the left end point and then increases to the right end point in $q$ and $1<\frac{2}{q(1-q)}\left(\frac{1-q}{2-q}\right)^{2-q}$. Hence, $y^{2-q}<\bar{\lambda}$ for all $q\in(0,1)$.

Since $h(x)$ has at most two stationary points, we denote them as $x{''} < x'$ if they exist. The function $h(x)$ increases from $x=0$ to $x=x{''}$ (local maximizer), then decreases until $x=x'$ (local minimizer), and finally increases until $x=y$.

To determine the global minimizer, we compare $h(0) = y^2$, $h(x')$ (if it exists), and $h(y) = \lambda y^q$:

\begin{enumerate}
  \item[(i)] If $\lambda < y^{2-q}$, then $h(y) = \lambda y^q < y^2 = h(0)$, and $h(x') < h(y)$, so $x'$ is the unique global minimizer.
  
  \item[(ii)] If $\lambda = y^{2-q}$, then $h(y) = \lambda y^q = y^2 = h(0)$ and $h(x') < h(0)$, which implies that $x'$ is the unique global minimizer.

  \item[(iii)] If $y^{2-q}<\lambda <\bar{\lambda}$, then $h(y) = \lambda y^q > y^2 = h(0)$:
  \begin{itemize}
    \item If  $h(x') < h(0)$, then $x'$ is the unique global minimizer.
    \item If $h(x') = h(0)$, then both $x=0$ and $x=x'$ are global minimizers.
    \item If $h(x') > h(0)$, then $x=0$ is the unique global minimizer.
  \end{itemize}
  \item[(iv)] If $\bar{\lambda} \leq\lambda$, then $x=0$ is the unique global minimizer.  
\end{enumerate}

In case (iii), we can further determine a sub-interval within $(y^{2-q},\bar{\lambda})$ to determine the global minimizer. The critical point $x'$ satisfies $h'(x')=0$. For $x\in (0,y)$, this gives
\begin{equation}\label{equ:lem1}
    \lambda = \frac{2(y-x')}{q(x')^{q-1}}\qquad \forall q \in (0,1). 
\end{equation}
We wish to find a threshold $\lambda^{\ast}\in(y^{2-q},\bar{\lambda})$ such that $h(x')=h(0)$ at $\lambda=\lambda^{\ast}$:
\begin{equation*}
    (y-x')^2+\lambda(x')^q =y^2,
\end{equation*}
which implies that
\begin{equation*}
    x'=\frac{2y(1-q)}{2-q}.
\end{equation*}
Substituting $x'$ back into \eqref{equ:lem1} gives
\[
\lambda^{\ast} =\frac{4}{2^q(1-q)}\cdot\left(\frac{1-q}{2-q} y\right)^{2-q} \qquad \forall q \in (0,1).
\]
It is obvious that $\lambda^{\ast}\in(y^{2-q},\bar{\lambda})$. Hence we obtain 
\begin{enumerate}
    \item[(iii')] If $y^{2-q}<\lambda\leq \lambda^{\ast} $, then $h(y) = \lambda y^q > y^2 = h(0)$:
  \begin{itemize}
    \item If  $h(x') < h(0)$, then $x'$ is the unique global minimizer.
    \item If $h(x') = h(0)$, then both $x=0$ and $x=x'$ are global minimizers.
  \end{itemize} 
  \item[(iv')] If $\lambda^{\ast}<\lambda$, then $x=0$ is the unique global minimizer.
\end{enumerate}

In addition, if $x'$ exists, the seconder order bound holds:
\[
h''(x')= 2+\lambda q(q-1)(x')^{q-2}> 0, \qquad \forall \, \lambda <\bar{\lambda},
\]
which implies that 
\[
x' > \left(\frac{\lambda q (1-q)}{2}\right)^{\frac{1}{2-q}}.
\]
\end{proof}
\begin{remark}
When $q$ equals $\frac{1}{2}$ or $\frac{1}{3}$,  closed-form expressions of $\psi_H(y;\lambda)$ are available in references \cite{Xu2013fast,Xu2012half}. Finding the minimizer of this nonconvex function is also investigated in \cite{peng2025newton}.
\end{remark}

\begin{remark}
Among various methods for finding roots (if they exist) of $h'(x)=0$, Newton's method is a good choice here because it provides rapid local convergence for polynomial functions.
\end{remark}

\begin{remark}
In the above proof, we say that 
\[\frac{2}{q(1-q)}\left( \frac{1-q}{2-q} \right)^{2-q}\]
decreases from the left end point and then increases to the right end point in $q\in(0,1)$ and $1<\frac{2}{q(1-q)}\left( \frac{1-q}{2-q} \right)^{2-q}$ for $q\in(0,1)$. Indeed, let 
\[
g(q)= \frac{2}{q(1-q)}\left( \frac{1-q}{2-q} \right)^{2-q} \qquad \forall q \in(0,1).
\]
We obtain
\[
\ln g(q) = \ln 2 - \ln q - \ln(1-q) + (2-q)[\ln(1-q)-\ln(2-q)],
\]
which implies that
\[
\frac{\operatorname{d}}{\operatorname{d}q}\ln g(q) = -\frac{1}{q} - \ln(1-q)+\ln(2-q)
\]
and 
\[
\frac{\operatorname{d}^2}{\operatorname{d}q^2}\ln g(q) = \frac{1}{q^2} + \frac{1}{1-q}-\frac{1}{2-q} >0.
\]
Then we can numerically find the root $q^{\ast}\approx0.691766$ via Newton's method. Hence, $g(q)\geq g(q^{\ast})\approx 1.4154$ for all $q \in (0,1)$.
\end{remark}

\section{Nonconvex regularization approximation}
Motivated by the thresholding recovery principle in nonconvex optimization, we propose a family of recovery thresholding hyperinterpolations. These methods adapt thresholding operators to the hyperinterpolation framework, enforcing sparsity while suppressing noise in the coefficient domain. Specifically, we define three variants: hard thresholding hyperinterpolation \cite{an2023hard}, springback hyperinterpolation \cite{an2022springback}, Newton hyperinterpolation.

In practice, the sampling data $\{f(\mathbf{x}_j)\}_{j=1}^{N}$ at nodes $\{\mathbf{x}_j\}_{j=1}^{N}$ are perturbed by noise $\epsilon$. We usually obtain the data $\{\mathbf{x}_j, f^{\epsilon}(\mathbf{x}_j)\}_{j=1}^{N}$ with $f^{\epsilon}(\mathbf{x}_j)=f(\mathbf{x}_j)+ \epsilon_j$, and actually solve the following three nonconvex regularized least squares approximation problems
\begin{enumerate}
    \item[$\bullet$] $\ell_0$-regularized 
    \begin{equation}\label{l0regulariztion}
\min\limits_{p_{\lambda}\in\Pi_n(\Omega)}~~\left\{\sum_{j=1}^{N} w_j[p_{\lambda}(\mathbf{x}_j)-f^{\epsilon}(\mathbf{x}_j)]^2
+ \lambda^2 \sum_{\ell=1}^{d_n} |\beta_{\ell}^{\lambda}|_0\right\},
    \end{equation}
    \item[$\bullet$] $\ell_q$-regularized 
    \begin{equation}\label{lqregulariztion}
    \min\limits_{p_{\lambda}\in\Pi_n(\Omega)}~~\left\{\sum_{j=1}^{N} w_j[p_{\lambda}(\mathbf{x}_j)-f^{\epsilon}(\mathbf{x}_j)]^2
+ \lambda \sum_{\ell=1}^{d_n} |\beta_{\ell}^{\lambda}|^{q}\right\},
    \end{equation}
    
    \item[$\bullet$] springback-regularized 
    \begin{equation}\label{spregulariztion}
\min\limits_{p_{\lambda}\in\Pi_n(\Omega)}~~\left\{ \frac{1}{2}\sum_{j=1}^{N} w_j[p_{\lambda}(\mathbf{x}_j)-f^{\epsilon}(\mathbf{x}_j)]^2
+ \lambda \sum_{\ell=1}^{d_n} \left( |\beta_{\ell}^{\lambda}|- \frac{\alpha}{2}|\beta_{\ell}^{\lambda}|^2 \right) \right\},
    \end{equation}
where $q\in(0,1)$, $ p_{\lambda}(\mathbf{x})=\sum_{\ell=1}^{d_n}\beta_{\ell}^{\lambda} \Phi_{\ell}(\mathbf{x}) \in \Pi_{n}(\Omega)$, $\lambda>0$ is the \emph{regularization parameter}, and $\alpha>0$ is a model parameter, and $|\beta_{\ell}^{\lambda}|$ represents the absolute value of $\beta_{\ell}^{\lambda}$, and $|\beta_{\ell}^{\lambda}|_{0}$ denotes the $\ell_0$-norm in one-dimension case, that is
\begin{equation*}
    \forall \beta_{\ell}^{\lambda} \in \mathbb{R},\quad |\beta_{\ell}^{\lambda}|_{0}:=
\left\{
\begin{array}{ll}
0, & \text{if}~~ \beta_{\ell}^{\lambda}=0, \\
1, & \text{if}~~ \beta_{\ell}^{\lambda}\ne0.
\end{array}
\right.
\end{equation*}
\end{enumerate}

We define \emph{hard thresholding hyperinterpolation}, \emph{Newton hyperinterpolation}, and \emph{springback hyperinterpolation} as follows:
\begin{definition}[Hard thresholding hyperinterpolation]
Let $\Omega \subset {\mathbb{R}}^s$  be a compact domain and $f \in \mathcal{C}(\Omega)$.
Suppose that the discrete scalar product $\langle f,g \rangle_N$ in \eqref{equ:discreteinner} is determined by an $N$-point quadrature rule of PI-type in $\Omega$ with algebraic degree of exactness $2n$. The \emph{hard thresholding hyperinterpolation} of $f$ onto $\Pi_{n}(\Omega)$ is defined as
\begin{eqnarray}\label{equ:hardhyper}
  \mathcal{H}_{n,0}^{\lambda}{f}:=\sum_{\ell=1}^{d_n} \eta_H (\langle{f, \Phi_{\ell}} \rangle_N; \lambda)\Phi_{\ell}, \qquad \lambda>0,
\end{eqnarray}
where $\eta_H$ is a \emph{hard thresholding operator} 
\begin{equation*}
    \eta_{H}(\langle{f, \Phi_{\ell}} \rangle_N;\lambda):=\left\{\begin{array}{cl}
\langle{f, \Phi_{\ell}} \rangle_N,  & \text{if}~~|\langle{f, \Phi_{\ell}} \rangle_N|> \lambda ,\\
0,  & \text{if}~~|\langle{f, \Phi_{\ell}} \rangle_N|\leq \lambda.
\end{array}\right.
\end{equation*}
\end{definition}

\begin{definition}[Newton hyperinterpolation] Let $\Omega \subset {\mathbb{R}}^s$  be a compact domain and $f \in \mathcal{C}(\Omega)$.
Suppose that the discrete scalar product $\langle f,g \rangle_N$ in \eqref{equ:discreteinner} is determined by an $N$-point quadrature rule of PI-type in $\Omega$ with algebraic degree of exactness $2n$. Given $q\in(0,1)$, the \emph{Newton hyperinterpolation} of $f$ onto $\Pi_{n}(\Omega)$ is defined as
\begin{eqnarray}\label{equ:newtonhyper}
  \mathcal{H}_{n,q}^{\lambda}{f}:=\sum_{\ell=1}^{d_n} \psi_{H} (\langle{f, \Phi_{\ell}} \rangle_N; \lambda)\Phi_{\ell}, \qquad \lambda>0,
\end{eqnarray}
where $\psi_H$ is defined by \eqref{equ:newton} in Lemma \ref{lem:newton}.
\end{definition}

\begin{definition}[Springback hyperinterpolation]
Let $\Omega \subset {\mathbb{R}}^s$  be a compact domain and $f \in \mathcal{C}(\Omega)$.
Suppose that the discrete scalar product $\langle f,g \rangle_N$ in \eqref{equ:discreteinner} is determined by an $N$-point quadrature rule of PI-type in $\Omega$ with algebraic degree of exactness $2n$. The \emph{springback hyperinterpolation} of $f$ onto $\Pi_{n}(\Omega)$ is defined as
\begin{eqnarray}\label{equ:springbackhyper}
  \mathcal{H}_{n}^{\lambda,\alpha}{f}:=\sum_{\ell=1}^{d_n} s_H (\langle{f, \Phi_{\ell}} \rangle_N; \lambda,\alpha)\Phi_{\ell}, \qquad \lambda>0,
\end{eqnarray}
where $s_H$ is a \emph{springback thresholding operator}  and $1-\lambda \alpha >0$
\begin{equation*}
    s_{H}(\langle{f, \Phi_{\ell}} \rangle_N;\lambda, \alpha):=\left\{\begin{array}{cl}
{\operatorname{sign}}(\langle{f, \Phi_{\ell}} \rangle_N) \frac{|\langle{f, \Phi_{\ell}} \rangle_N|- \lambda}{1-\lambda \alpha},  & \text{if}~~|\langle{f, \Phi_{\ell}} \rangle_N|> \lambda ,\\
0,  & \text{if}~~|\langle{f, \Phi_{\ell}} \rangle_N|\leq \lambda. 
\end{array}\right.
\end{equation*}
\end{definition}

For comparison with these nonconvex approximations, we include \emph{Lasso hyperinterpolation}, a convex approximation method that provides the unique solution to an $\ell_1$-regularized least squares problem \cite{an2021lasso}.

\begin{definition}[Lasso hyperinterpolation]
Let $\Omega \subset {\mathbb{R}}^s$  be a compact domain and $f \in \mathcal{C}(\Omega)$.
Suppose that the discrete scalar product $\langle f,g \rangle_N$ in \eqref{equ:discreteinner} is determined by an $N$-point quadrature rule of PI-type in $\Omega$ with algebraic degree of exactness $2n$. The \emph{Lasso hyperinterpolation} of $f$ onto $\Pi_{n}(\Omega)$ is defined as
\begin{eqnarray}\label{equ:lassohyper}
  \mathcal{H}_{n,1}^{\lambda}{f}:=\sum_{\ell=1}^{d_n} \eta_S (\langle{f, \Phi_{\ell}} \rangle_N; \lambda )\Phi_{\ell}, \qquad \lambda>0,
\end{eqnarray}
where $\eta_S$ is a \emph{soft thresholding operator} defined as 
\begin{equation*}
    \eta_{S}(\langle{f, \Phi_{\ell}} \rangle_N;\lambda):=\left\{\begin{array}{cl}
{\operatorname{sign}}(\langle{f, \Phi_{\ell}} \rangle_N) (|\langle{f, \Phi_{\ell}} \rangle_N|- \lambda),  & \text{if}~~|\langle{f, \Phi_{\ell}} \rangle_N|> \lambda ,\\
0,  & \text{if}~~|\langle{f, \Phi_{\ell}} \rangle_N|\leq \lambda. 
\end{array}\right.
\end{equation*}
\end{definition}

\begin{remark}
For latter use, we list different notations to distinguish various hyperinterpolation:
\begin{table}[H]
\begin{center}
  \begin{tabular}{|c|c|c|c|c|}
\hline 
     &   \text{Hard thresholding}& \text{Newton} $q \in (0,1) $  & \text{Lasso} &\text{Springback}   \\ \hline
 \textbf{Classical}   & $\mathcal{H}_{n,0}^{\lambda} f$&$\mathcal{H}_{n,q}^{\lambda} f$ &$\mathcal{H}_{n,1}^{\lambda} f$ &$\mathcal{H}_{n}^{\lambda,\alpha} f$  \\
\textbf{Relaxed}   & $\mathfrak{H}_{n,0}^{\lambda} f$&$\mathfrak{H}_{n,q}^{\lambda} f$ &$\mathfrak{H}_{n,1}^{\lambda}f$&$\mathfrak{H}_{n}^{\lambda,\alpha} f$   \\
\textbf{Unfettered}   & $\mathrm{H}_{n,0}^{\lambda} f$& $\mathrm{H}_{n,q}^{\lambda} f$ &$\mathrm{H}_{n,1}^{\lambda}f$ &$\mathrm{H}_{n}^{\lambda,\alpha} f$  \\

 \hline 
  \end{tabular}
  \caption{Summary of notations for various hyperinterpolation methods.}
  \label{tab:various_hyper}
\end{center}
\end{table}
\end{remark}

Now, we obtain the following important result: 
\begin{theorem}\label{thm:ExactSolution}
Let $\Omega \subset {\mathbb{R}}^s$ be a compact domain and $f \in \mathcal{C}(\Omega)$. Let $\alpha >0$ be such that $1-\lambda \alpha >0$.
Suppose that the discrete scalar product $\langle f,g \rangle_N$ in \eqref{equ:discreteinner} is determined by an $N$-point quadrature rule of PI-type in $\Omega$ with algebraic degree of exactness $2n$.   Then $\mathcal{H}_{n,0}^{\lambda} f$,  $\mathcal{H}_{n,q}^{\lambda}f$ and $\mathcal{H}_{n}^{\lambda,\alpha}f$  are the unique solutions to the nonconvex regularized least squares approximation problems \eqref{l0regulariztion}, \eqref{lqregulariztion} with $q\in(0,1)$ and \eqref{spregulariztion} with noise $\epsilon \equiv 0$, respectively.
\end{theorem}

\begin{proof}
Let $\alpha_{\ell}=\langle f , \Phi_{\ell}\rangle_N$ for $\ell=1,\cdots, d_n$. Since the case of $\ell_0$-regularized approximation has been proved in \cite{an2023hard}, we only need solve the following two one-dimensional cases:
\[
   h(\beta_{\ell}^{\lambda}):=( \beta_{\ell}^{\lambda} - \alpha_{\ell})^2 + \lambda |\beta_{\ell}^{\lambda}|^q \qquad \forall q\in (0,1),
\]
and 
\[
  g(\beta_{\ell}^{\lambda}) :=\frac{1}{2} ( \beta_{\ell}^{\lambda} - \alpha_{\ell})^2 + \lambda \left( |\beta_{\ell}^{\lambda}| - \frac{\alpha}{2} |\beta_{\ell}^{\lambda}|^2  \right).
\]
We can find the minimizer of $h(\beta_{\ell}^{\lambda})$ by \eqref{equ:lem1} in Lemma \ref{lem:newton}. 

Now we focus on finding the minimizer of $ g(\beta_{\ell}^{\lambda})$. Taking the first derivative of $g(\beta_{\ell}^{\lambda})$ and setting it equal to zero yields
\[
g'(\beta_{\ell}^{\lambda}) = \beta_{\ell}^{\lambda} - \alpha_{\ell} + \lambda\partial(|\beta_{\ell}^{\lambda}|) - {\lambda \alpha} \partial(|\beta_{\ell}^{\lambda}|) |\beta_{\ell}^{\lambda}|=0,
\]
where 
\[
\partial(|\beta_\ell^{\lambda}|) = \begin{cases}
1 & \text{if } \beta_\ell^{\lambda} > 0, \\
-1 & \text{if } \beta_\ell^{\lambda} < 0, \\
\in [-1, 1] & \text{if } \beta_\ell^{\lambda} = 0.
\end{cases}
\]
Thus, we have 
\[
\beta^{\lambda}_{\ell} =  \begin{cases}
\operatorname{sign}(\alpha_{\ell}) \frac{|\alpha_{\ell}| - \lambda}{1-\lambda \alpha} & \text{if } |\alpha_\ell| > \lambda, \\
0 & \text{if } |\alpha_\ell| \leq \lambda.
\end{cases}
\]
Therefore, we have completed the proof.
\end{proof}

\section{Sparsity analysis}
Motivated by Theorem 1 in \cite{Chen2014complexity_nonconvex}, we analyze the sufficient condition on $\lambda$ for global minimizers of the nonconvex regularization approximation to have desirable sparsity.
\begin{theorem}
Let $F_q(\bm{\beta}^{\lambda} ) = \|\mathbf{W}^{\frac{1}{2}}(\mathbf{A}\bm{\beta}^{\lambda} - \mathbf{f}^{\epsilon}  )\|_2^2 + \lambda \| \bm{\beta}^{\lambda}\|_q^q $ for $q\in (0,1)$, and let $\bm{\alpha}=[\alpha_1, \cdots, \alpha_{d_n}]^{\operatorname{T}}\in \mathbb{R}^{d_n}$ with $\alpha_{\ell}= \langle f^{\epsilon}, \Phi_{\ell} \rangle_N$ for $\ell=1,\cdots,d_n$. Let $\mathbf{u}=[u_1,\cdots,u_{d_n}]^{\operatorname{T}}\in \mathbb{R}^{d_n}$ with
\[
u_\ell=\frac{4}{2^q(1-q)} \left(\frac{1-q}{2-q} |\alpha_{\ell}|\right)^{2-q} \qquad \forall q\in (0,1).
\]
Let $u_{(1)},u_{(2)}, \cdots,u_{(d_n)}$ denote the components of $\mathbf{u}$ rearranged in descending order such that
\[
u_{(1)}\geq u_{(2)} \geq  \cdots \geq u_{(d_n)}.
\]
Given an integer $k\geq 2$, if $ \lambda \geq u_{(k) }$,  any global minimizer $\bm{\beta}^{\lambda}=[\beta_1^{\lambda},\cdots,\beta_{d_n}^{\lambda}]^{\operatorname{T}} \in \mathbb{R}^{d_n}$ of problem \eqref{lqregulariztion}   satisfies $\|\bm{\beta}^{\lambda}\|_0 \leq k-1$.
\end{theorem}

\begin{proof}
Since the objective function $F_q(\bm{\beta}^{\lambda} )$ is separable across coordinates, the global minimizer is given by the component-wise minimizer of $(\beta_{\ell}^{\lambda} - \alpha_{\ell})^2 + \lambda|\beta_{\ell}^{\lambda}|^{q}$ for each $\ell$:
\[
\begin{aligned}
    F_q(\bm{\beta}^{\lambda} ) = \|\mathbf{W}^{\frac{1}{2}}(\mathbf{A}\bm{\beta}^{\lambda} - \mathbf{f}^{\epsilon}  )\|_2^2 + \lambda \| \bm{\beta}^{\lambda}\|_q^q 
= C+\sum_{\ell=1}^{d_n}[(\beta_{\ell}^{\lambda} - \alpha_{\ell})^2 + \lambda|\beta_{\ell}^{\lambda}|^{q}],
\end{aligned}
\]
where $C=\sum_{j=1}^{N}w_jf^{\epsilon}[(\mathbf{x}_j)]^2 - \sum_{\ell=1}^{d_n}\alpha_{\ell}^2$. Then, by Lemma \ref{lem:newton}, the result is obvious.
\end{proof}

In particular, to obtain the desired sparsity $\|\bm{\beta}^{\lambda}\|_0=k$, we have
\[
\mathbb{P}(|\beta_\ell^{\lambda}|_0=1)=\frac{k}{d_n}\quad \text{and} \quad \mathbb{P}(|\beta_\ell^{\lambda}|_0=0)=1-\frac{k}{d_n}.
\]
Then by Bernstein’s inequality for bounded distributions \cite{Vershynin2018prob}, we have 
\[
\mathbb{P}(|\|\bm{\beta}^{\lambda}\|_0-k|> t)\leq 2\exp(-\frac{t^2/2}{\sigma^2+t/3}),
\]
where $\|\bm{\beta}^{\lambda}\|_0=\sum_{\ell=1}^{d_n}|\beta_{\ell}^{\lambda}|_0$ and $\sigma^2=\sum_{\ell=1}^{d_n}\operatorname{Var}(|\beta_{\ell}^{\lambda}|_0)=\operatorname{Var}(\|\bm{\beta}^{\lambda}\|_0)$.

\begin{lemma}\label{lem:4-1}
Let $\lambda > 0$ and $\eta$ be a sub-Gaussian random variable with $\mathbb{E}[\eta] = 0$ and $\|\eta\|_{\psi_2} < \infty$, where 
$\|\eta\|_{\psi_2}$ is the \emph{sub-Gaussian norm} of a random variable $\eta$, defined by 
\[\|\eta\|_{\psi_2}:= \inf \{t >0: \, \mathbb{E} \exp(\eta^2/t^2) \leq 2\}.\]
Define $\tilde{\alpha} = \alpha + \eta$. We have the following results.
\begin{enumerate}
    \item If $|\alpha| > \lambda$, then 
    \begin{equation*}
        \mathbb{P}(|\tilde{\alpha}| < \lambda) \leq 2\exp\left(-\frac{(|\alpha| - \lambda)^2}{2\|\eta\|_{\psi_2}^2}\right).
    \end{equation*}

    \item If $|\alpha| \leq \lambda$, then 
    \begin{equation*}
        \mathbb{P}(|\tilde{\alpha}| > \lambda) \leq 2\exp\left(-\frac{(\lambda - |\alpha|)^2}{2\|\eta\|_{\psi_2}^2}\right).
    \end{equation*}
\end{enumerate}
\end{lemma}

\begin{proof}
We use the standard tail bound for sub-Gaussian variables: if $\eta$ has $\mathbb{E}[\eta] = 0$ and sub-Gaussian norm $\|\eta\|_{\psi_2}<\infty$, then
\begin{equation*}
    \mathbb{P}(|\eta| \geq t) \leq 2\exp\left(-\frac{t^2}{2\|\eta\|_{\psi_2}^2}\right)\quad\forall t \geq 0.
\end{equation*}

\textbf{Case 1:} Suppose $|\alpha| > \lambda$.
Note that  $\{|\alpha+\eta|<\lambda\} \subset \{|\eta| > |\alpha|-\lambda\}$ implies that
\begin{equation*}
    \mathbb{P}(|\alpha + \eta| < \lambda) \leq \mathbb{P}(|\eta| > |\alpha|-\lambda )\leq 2\exp\left(-\frac{(|\alpha| - \lambda)^2}{2\|\eta\|_{\psi_2}^2}\right).
\end{equation*}

\textbf{Case 2:} Suppose $|\alpha| \leq \lambda$. The event ${|\alpha + \eta| > \lambda}$ occurs when either $\alpha + \eta > \lambda$ or $\alpha + \eta < -\lambda$. This gives
\begin{equation*}
    \mathbb{P}(|\alpha + \eta| > \lambda) \leq \mathbb{P}(\eta > \lambda - \alpha) + \mathbb{P}(\eta < -\lambda - \alpha).
\end{equation*}
For the first term, since $\lambda - \alpha \geq \lambda-|\alpha|\geq 0$, we have 
\[
   \mathbb{P}(\eta>\lambda - \alpha) \leq    \mathbb{P}(\eta > \lambda - |\alpha|),% \leq    \mathbb{P}(|\eta|\geq \lambda - |\alpha|),
\]
and similarly for the second term by $-\lambda - \alpha \leq -(\lambda - |\alpha|))$,
\begin{equation*}
        \mathbb{P}(\eta < -\lambda - \alpha) \leq \mathbb{P}(\eta < -(\lambda - |\alpha|)) .% \leq \mathbb{P}(|\eta| \geq \lambda - |\alpha|).
\end{equation*}

Thus combining both terms by one-sided bounds for sub-Gaussian random variables:
\begin{equation*}
\mathbb{P}(|\alpha + \eta| > \lambda) \leq \mathbb{P}(\eta \geq \lambda - |\alpha|) + \mathbb{P}(\eta < -(\lambda - |\alpha|))\leq 2\exp\left(-\frac{(\lambda - |\alpha|)^2}{2\|\eta\|_{\psi_2}^2}\right).
\end{equation*}
\end{proof}

\begin{corollary}\label{cor:4-1}
Under conditions of Lemma \ref{lem:4-1}, for $\delta \in (0,1)$ and $c > 0$, if
\[
\left| |\alpha| - \lambda \right| > c\lambda \quad \text{and} \quad \lambda \geq \frac{\|\eta\|_{\psi_2}}{c} \sqrt{2\ln \frac{2}{\delta}},
\]
then the probability of threshold decision mismatch satisfies:
\[
\mathbb{P} \left[ (|\alpha|-\lambda)(|\tilde{\alpha}-\lambda|) <0 \right] \leq \delta.
\]
\end{corollary}

\begin{proof}
By Lemma \ref{lem:4-1} and the condition $\left| |\alpha| - \lambda \right| > c\lambda$:
\[
 \mathbb{P}\left[ (|\alpha|-\lambda)(|\tilde{\alpha}-\lambda|) <0 \right] \leq 2\exp\left(-\frac{\left| |\alpha| - \lambda \right|^2}{2\|\eta\|_{\psi_2}^2}\right) \leq 2\exp\left(-\frac{(c\lambda)^2}{2\|\eta\|_{\psi_2}^2}\right).
\]
Substituting $\lambda \geq \frac{\|\eta\|_{\psi_2}}{c} \sqrt{2\ln \frac{2}{\delta}}$:
\[
2\exp\left(-\frac{c^2 \lambda^2}{2\|\eta\|_{\psi_2}^2}\right) \leq 2\exp\left(-\ln \frac{2}{\delta}\right) = \delta.
\]
Thus, we have completed the proof.
\end{proof}

%%%%%%%%%%%%%%%%%%%%%%%%%%%%%%%%%%%%%%%%%%%%%%%%%

In the following, we only consider hard thresholding hyperinterpolation. Let $\tilde{\bm{\alpha}}= \bm{\alpha}+\bm{\eta}$, where $\bm{\alpha}=\mathbf{A}^{\rm{T}}\mathbf{Wf}$ and $\bm{\eta}=\mathbf{A}^{\rm{T}}\mathbf{W}\bm{\epsilon}  $. We assume that $\bm{\epsilon}=[\epsilon_1,\cdots, \epsilon_N]^{\rm{T}}$ is mean zero and sub-Gaussian, i.e., $\mathbb{E}[\epsilon_j]=0$ and $\|\epsilon_j\|_{\psi_2}<\infty$ for $j=1,\cdots,N$. Let $R = \max_{1\leq j \leq N}\|\epsilon_j\|_{\psi_2}$. We know that $\mathbb{E}[\eta_{\ell}]=0$ and
\[
\|\eta_{\ell}\|_{\psi_2}^2=\left\| \sum_{j=1}^{N}w_j \Phi_{\ell}(\mathbf{x}_j)\epsilon_j \right\|_{\psi_2}^2 \leq CR^2 \|\mathbf{A}^{\rm{T}} \mathbf{W}\|_2^2<\infty  , \qquad \forall \ell =1, \cdots, d_n.
\]

Let $S_K= \mathrm{supp}(\bm{\alpha}_K)$ be the oracle support set of size $K$ as the magnitude of the $K$-th largest coefficient in ${\bm{\alpha}}$, and let $\hat{S}_K$ be the support of $\tilde{\bm{\alpha}}$ selected by hard thresholding operators. Define $\lambda_K$ as the magnitude of the $K$-th largest coefficient in $\tilde{\bm{\alpha}}$.

%%%%%%%%%%%%%%%%%%%%%%%%%%%%%%%%%%%%%%%%%%%%%%%%%

\begin{theorem}[Approximation Error Bound for Sub-Gaussian Noise]\label{thm:4-2}
Let $\bm{\alpha} \in \mathbb{R}^{d_n}$ be the true coefficient vector and $\tilde{\bm{\alpha}} = \bm{\alpha}+\bm{\eta}$ be the noisy observation, where $\bm{\eta}=[\eta_1,\cdots,\eta_{d_n}]^{\rm{T}} \in \mathbb{R}^{d_n}$, and $\{\eta_\ell\}_{\ell=1}^{d_n}$ are independent sub-Gaussian random variables with $\mathbb{E}[\eta_{\ell}] = 0$ and $\|\eta_{\ell}\|_{\psi_2} < \infty$. Let $R=\max_{1\leq \ell \leq d_n} \|\eta_{\ell}\|_{\psi_2}$. Let $\bm{\beta}^{\lambda}_K$ be the coefficients retained by hard thresholding operator that retains the $K$ largest coefficients in magnitude, i.e., $\beta_{\ell}^{\lambda}=\tilde{\alpha}_{\ell}$ if $|\tilde{\alpha}_{\ell}|>\lambda$. Then 
\begin{equation}
\mathbb{E}[\|\bm{\beta}^{\lambda}_K -\bm{\alpha}\|_2^2] \leq C KR + \| \bm{\alpha}_{\hat{S}_K^c}\|_2^2,
\end{equation}
where $C>0$ is a universal constant.
\end{theorem}

\begin{proof}
Note that
\begin{equation}\label{equ:4-1}
\|\bm{\beta}^{\lambda}_K -\bm{\alpha}\|_2^2 = \sum_{\ell \in \hat{S}_K} |\eta_{\ell}|^2 + \sum_{\ell \notin \hat{S}_K} |\alpha_{\ell}|^2.
\end{equation}
Since there exists a constant $C$ such that
\[
\mathbb{E}[|\eta_{\ell}|^2] \leq C \|\eta_{\ell}\|_{\psi_2},
\]
we have 
\[
\mathbb{E}\left[\sum_{\ell \in \hat{S}_K} |\eta_{\ell}|^2 \right] \leq CKR.
\]
For the second sum in \eqref{equ:4-1}, we have 
\[
\| \bm{\alpha}_{\hat{S}_K^c}\|_2^2=\sum_{\ell \in  \hat{S}_K^c} |\alpha_{\ell}|^2 .
\]
Thus, we have completed the proof.
\end{proof}

\begin{theorem}\label{thm:4-3}
Under the same conditions of Theorem \ref{thm:4-2}, let $\tilde{\bm{\gamma}}_{q,K}^{\lambda_1}=\psi_H(\tilde{\bm{\alpha}};\lambda_1)$ be the $K$ largest coefficients in magnitude  defined by \eqref{equ:newton} in Lemma \ref{lem:newton}. Then for Newton hyperinterpolation, we have 
\begin{equation}\label{equ:4-2}
\mathbb{E}[ \| \mathcal{H}_{n,q}^{\lambda_1} f^{\epsilon} -f\|_2] \leq  \mathbb{E} [ \|\tilde{\bm{\gamma}}_{q,K}^{\lambda_1} - {\bm{\beta}}_{K}^{\lambda_2}\|_2]  + \sqrt{C KR + \| \bm{\alpha}_{\hat{S}_K^c}\|_2^2} + 2V^{\frac{1}{2}}E_n(f).
\end{equation} 
where $E_n(f):=\inf\{\|f-p\|_{\infty}; p \in \mathbb{P}_n(\Omega)\}$, and $V=\int_{\Omega} \, \text{d} \omega$, and $C>0$ is a universal constant.
\end{theorem}

\begin{proof}
First, we have the following decomposition:
\[
 \| \mathcal{H}_{n,q}^{\lambda_1} f^{\epsilon} -f\|_2 \leq  \| \mathcal{H}_{n,q}^{\lambda_1} f^{\epsilon} - \mathcal{H}_{n,0}^{\lambda_2}f\|_2 +  \| \mathcal{H}_{n,0}^{\lambda_2}f - \mathcal{L}_n f\|_2 + \| \mathcal{L}_n f -f\|_2,
\]
where $\| \mathcal{L}_n f -f\|_2 \leq 2V^{\frac{1}{2}} E_n(f)$ follows from \cite[Theorem 1]{sloan1995hyperinterpolation}.
By Parseval's identity, we have 
\[
 \| \mathcal{H}_{n,q}^{\lambda_1} f^{\epsilon} - \mathcal{H}_{n,q}^{\lambda_2}f\|_2 = \|\tilde{\bm{\gamma}}_{q,K}^{\lambda_1} - {\bm{\beta}}_{K}^{\lambda_2}\|_2 \quad \text{and} \quad  \| \mathcal{H}_{n,q}^{\lambda_2}f - \mathcal{L}_n f\|_2 = \| {\bm{\beta}}_{K}^{\lambda_2} - \bm{\alpha}\|_2.
\]
By Theorem \ref{thm:4-2} and Jensen's inequality, we have
\[
 \mathbb{E}[ \| {\bm{\beta}}_{K}^{\lambda_2} - \bm{\alpha}\|_2]  \leq \sqrt{ \mathbb{E}[ \| {\bm{\beta}}_{K}^{\lambda_2} - \bm{\alpha}\|_2^2]  } \leq \sqrt{C KR + \| \bm{\alpha}_{\hat{S}_K^c}\|_2^2}.
\]
Thus, we have completed the proof.
\end{proof}

\section{Numerical examples}
In particular, we considered
\begin{itemize}
\item[$\bullet$] {\em{Gaussian noise}} ${\cal{N}}(0,\sigma^2)$ from a normal distribution with mean 0
        and standard deviation {\tt{sigma}}=$\sigma$, implemented via the Matlab command
\begin{center}
{\tt{sigma*randn(N,1)}}.
\end{center}
\item[$\bullet$]  {\em{Impulse noise}} ${\cal{I}}(a)$ that takes a uniformly distributed random values in $[-a,a]$ with probability density $1/(2a)$ by means of the Matlab command
\begin{center}
{\tt{a*(1-2*rand(N,1)).*binornd(1,0.5,N,1)}},
\end{center}
where {\tt{binornd(1,0.5,N,1)}} generates an array of $N \times 1$ random binary numbers (0 or 1), with each number having the probability $1/2$ of being 1 and the probability $1/2$ of being 0.
\end{itemize}

\subsection{Signal processing}
\begin{figure}[htbp]
  \centering
  \includegraphics[width=\textwidth]{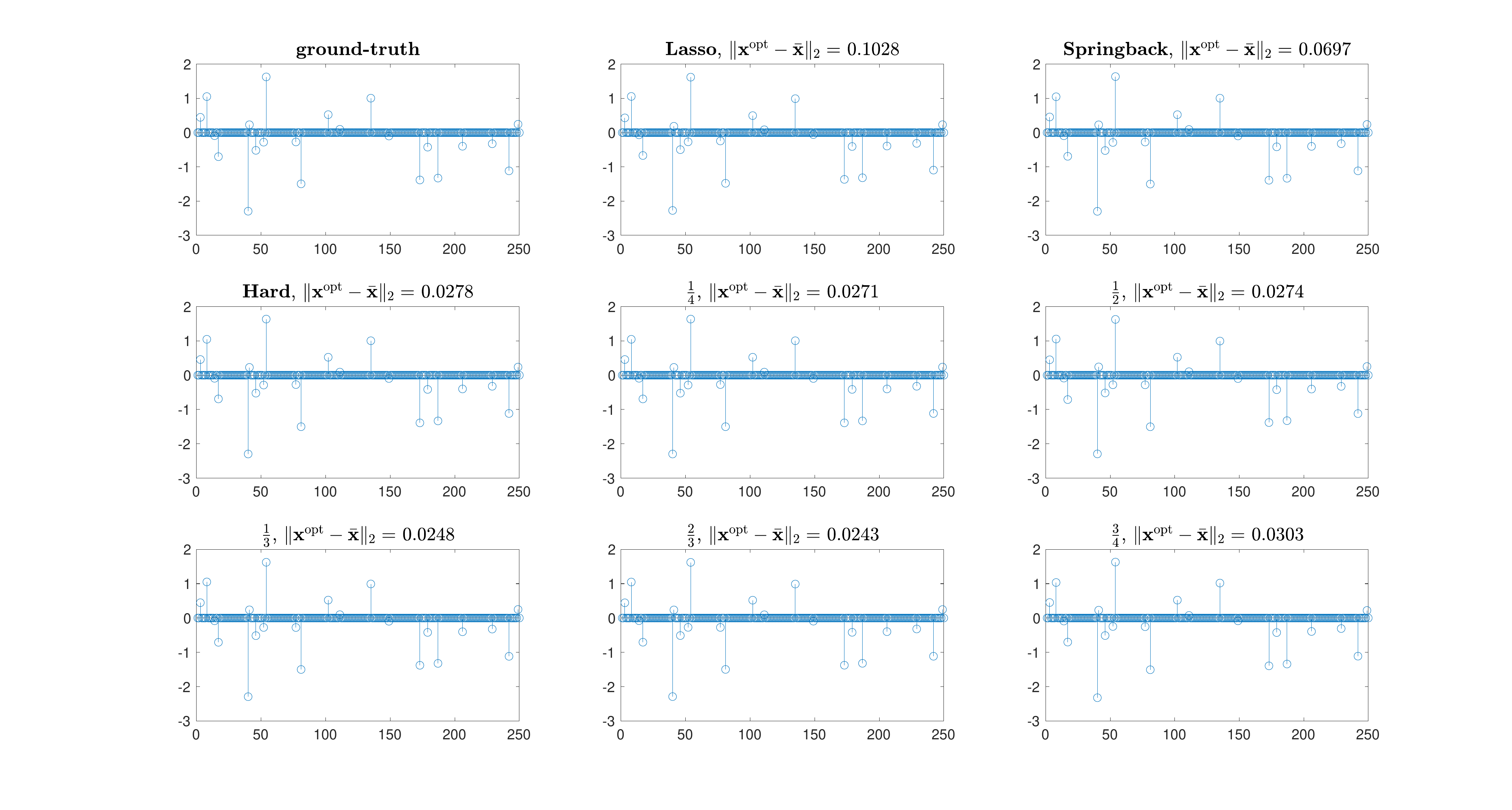}\\
  \caption{A ground-truth and its reconstructions from noisy measurements using Lasso hyperinterpolation, springback hyperinterpolation with $\alpha =1$, hard thresholding hyperinterpolation, Newton hyperinterpolation with $q = \frac{1}{4},\frac{1}{3},\frac{1}{2},\frac{2}{3},\frac{3}{4}$.}\label{fig:signal_processing}
\end{figure}

Let the sensing matrix $\mathbf{A}\in \mathbb{R}^{N\times d_n}$ with $(N,d_n)=(301,250)$, where $[\mathbf{A}]_{j\ell}=\Phi_{\ell}(\mathbf{x}_j)$ , 0the ground-truth $\bar{\mathbf{x}} \in \mathbb{R}^{250}$ be a 22-sparse vector with nonzero entries drawn from the standard normal distribution. The measurement vector $\mathbf{b}=\mathbf{A}\bar{\mathbf{x}}$ is contaminated by Gaussian noise  $\mathcal{N}(0,0.15^2)$. The ground-truth and its reconstructions are displayed in Figure \ref{fig:signal_processing}.

\begin{table}[H]
\begin{center}
  \begin{tabular}{c|c c c }
\hline 

$\sigma=0.15$   & $\|\mathbf{x}^{\text{opt}} - \bar{\mathbf{x}} \|_2$ & $\|\mathbf{x}^{\text{opt}} - \bar{\mathbf{x}} \|_{\infty}$ & $\text{AISNR}$ \\ \hline 

$\text{Lasso}$      &                $0.180969$      &      $0.059998$  &        $-2.753909$  \\[2pt] 

$\text{Springback}$ &                $0.129684$      &      $0.055627$  &        $1.340419$ \\[2pt] 

$\text{Hard}$       &                $0.056672$      &      $0.026948$  &        $7.589600$ \\[2pt] 

$\frac{1}{4}$       &                $0.057261$      &      $0.027241$  &        $7.497753$ \\[2pt] 

$\frac{1}{3}$       &                $0.057641$      &      $0.027438$  &        $7.438865$ \\[2pt] 

$\frac{1}{2}$       &                $0.059500$      &      $0.028360$  &        $7.169396$ \\[2pt] 

$\frac{2}{3}$       &                $0.065859$      &      $0.031126$  &        $6.299782$ \\[2pt] 

$\frac{3}{4}$       &                $0.073627$      &      $0.033966$  &        $5.316021$ \\[2pt] 
 \hline 
  \end{tabular}
  \caption{The average of $L_2$ and maximum errors, and improved SNRs (AISNR), where the original average SNR is $30.435188$ for $20000$ tests, and the parameter is $\alpha=1$ in Springback hyperinterpolation.}
  \label{tab:sigma=0.15}
\end{center}
\end{table}

\begin{table}[H]
\begin{center}
  \begin{tabular}{c|c c c }
\hline 
$\sigma=0.20$   & $\|\mathbf{x}^{\text{opt}} - \bar{\mathbf{x}} \|_2$ & $\|\mathbf{x}^{\text{opt}} - \bar{\mathbf{x}} \|_{\infty}$ & $\text{AISNR}$ \\ \hline 

$\text{Lasso}$      &                $0.241938$      &      $0.080156$  &        $-3.012405$  \\[2pt] 

$\text{Springback}$ &                $0.161096$      &      $0.066630$  &        $1.045825$ \\[2pt] 

$\text{Hard}$       &                $0.075644$      &      $0.035925$  &        $7.442771$ \\[2pt] 

$\frac{1}{4}$       &                $0.076442$      &      $0.036347$  &        $7.335256$ \\[2pt] 

$\frac{1}{3}$       &                $0.077068$      &      $0.036674$  &        $7.253117$ \\[2pt] 

$\frac{1}{2}$       &                $0.080089$      &      $0.038173$  &        $6.877130$ \\[2pt] 

$\frac{2}{3}$       &                $0.089655$      &      $0.042266$  &        $5.812357$ \\[2pt] 

$\frac{3}{4}$       &                $0.100776$      &      $0.046201$  &        $4.738579$ \\[2pt] 
 \hline 
  \end{tabular}
  \caption{The average of $L_2$ and maximum errors, and improved SNRs (AISNR), where the original average SNR is $27.266262$ for $20000$ tests, and the parameter is $\alpha=1$ in Springback hyperinterpolation.}
  \label{tab:sigma=0.20}
\end{center}
\end{table}

\begin{table}[H]
\begin{center}
  \begin{tabular}{c|c c c }
\hline 

$\sigma=0.25$   & $\|\mathbf{x}^{\text{opt}} - \bar{\mathbf{x}} \|_2$ & $\|\mathbf{x}^{\text{opt}} - \bar{\mathbf{x}} \|_{\infty}$ & $\text{AISNR}$ \\ \hline 

$\text{Lasso}$      &                $0.268514$      &      $0.093331$  &        $-2.632364$  \\[2pt] 

$\text{Springback}$ &                $0.158577$      &      $0.072293$  &        $3.391672$ \\[2pt] 

$\text{Hard}$       &                $0.111282$      &      $0.061713$  &        $6.319284$ \\[2pt] 

$\frac{1}{4}$       &                $0.107852$      &      $0.055568$  &        $6.591239$ \\[2pt] 

$\frac{1}{3}$       &                $0.106757$      &      $0.053440$  &        $6.672516$ \\[2pt] 

$\frac{1}{2}$       &                $0.105529$      &      $0.049920$  &        $6.708596$ \\[2pt] 

$\frac{2}{3}$       &                $0.109507$      &      $0.050184$  &        $6.149091$ \\[2pt] 

$\frac{3}{4}$       &                $0.117750$      &      $0.053483$  &        $5.281246$ \\[2pt] 
 \hline 
  \end{tabular}
  \caption{The average of $L_2$ and maximum errors, and improved SNRs (AISNR), where the original average SNR is $26.733392$ for $20000$ tests, and the parameter is $\alpha=1$ in Springback hyperinterpolation.}
  \label{tab:sigma=0.25}
\end{center}
\end{table}

\begin{table}[H]
\begin{center}
  \begin{tabular}{c|c c c }
\hline 

$\sigma=0.30$   & $\|\mathbf{x}^{\text{opt}} - \bar{\mathbf{x}} \|_2$ & $\|\mathbf{x}^{\text{opt}} - \bar{\mathbf{x}} \|_{\infty}$ & $\text{AISNR}$ \\ \hline 

$\text{Lasso}$      &                $0.320870$      &      $0.109237$  &        $-1.273842$  \\[2pt] 

$\text{Springback}$ &                $0.246952$      &      $0.097109$  &        $1.714914$ \\[2pt] 

$\text{Hard}$       &                $0.161049$      &      $0.080760$  &        $4.613875$ \\[2pt] 

$\frac{1}{4}$       &                $0.157028$      &      $0.078915$  &        $4.846475$ \\[2pt] 

$\frac{1}{3}$       &                $0.156221$      &      $0.078576$  &        $4.905019$ \\[2pt] 

$\frac{1}{2}$       &                $0.156858$      &      $0.078560$  &        $4.917246$ \\[2pt] 

$\frac{2}{3}$       &                $0.165097$      &      $0.079891$  &        $4.557312$ \\[2pt] 

$\frac{3}{4}$       &                $0.175977$      &      $0.081550$  &        $4.050387$ \\[2pt] 
 \hline 
  \end{tabular}
  \caption{The average of $L_2$ and maximum errors, and improved SNRs (AISNR), where the original average SNR is $21.881685$ for $20000$ tests, and the parameter is $\alpha=1$ in Springback hyperinterpolation.}
  \label{tab:sigma=0.3}
\end{center}
\end{table}

\subsection{Function denoising}
\begin{figure}[htbp]
  \centering
  \includegraphics[width=\textwidth]{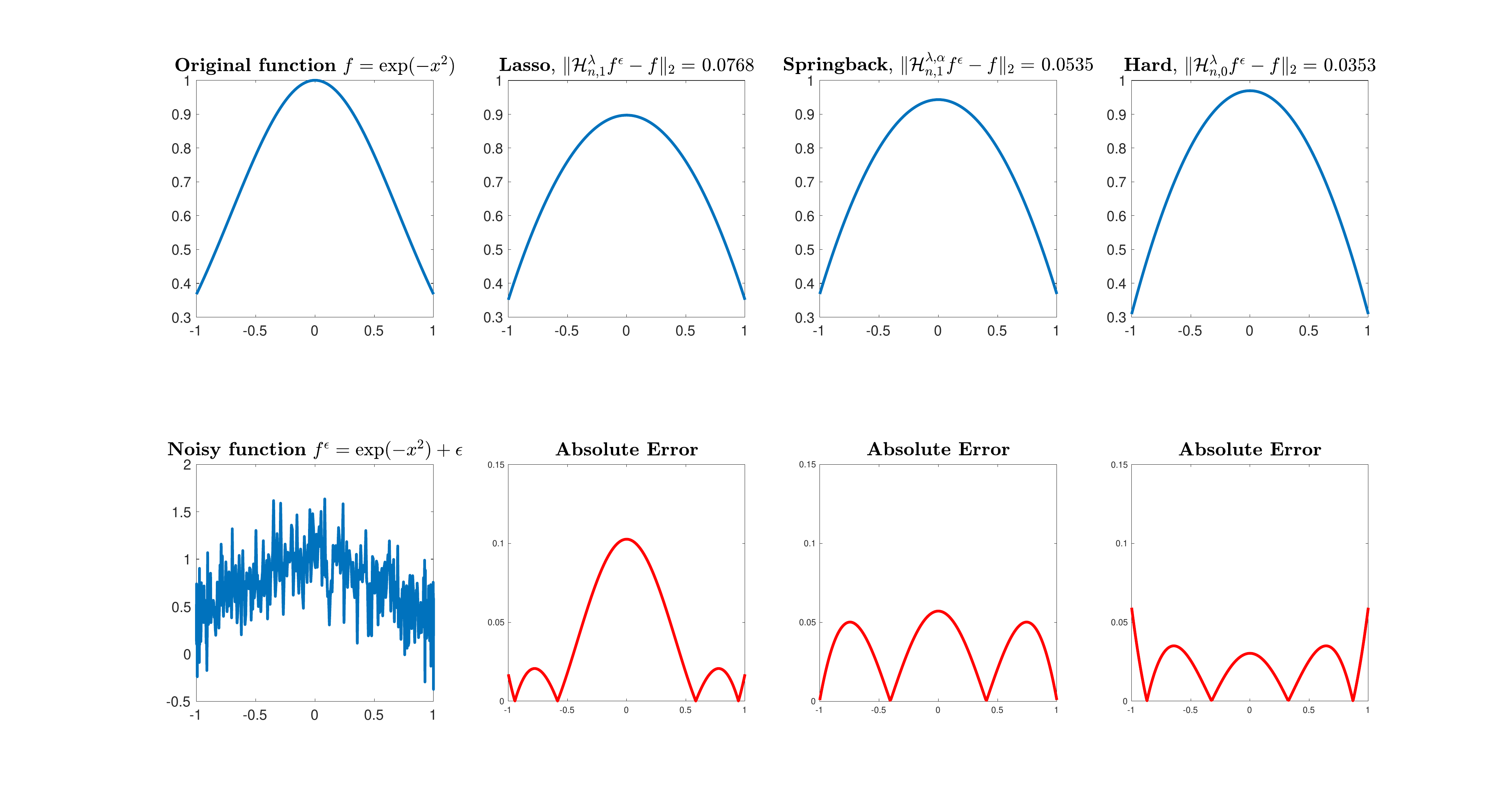}\\
  \caption{The original $f=\exp(-x^2)$ and its noisy version $f^{\epsilon}$, we reconstruct $f$ from noisy measurements using Lasso, springback, hard thresholding hyperinterpolations.}\label{fig:interval}
\end{figure}

\begin{figure}[htbp]
  \centering
  \includegraphics[width=\textwidth]{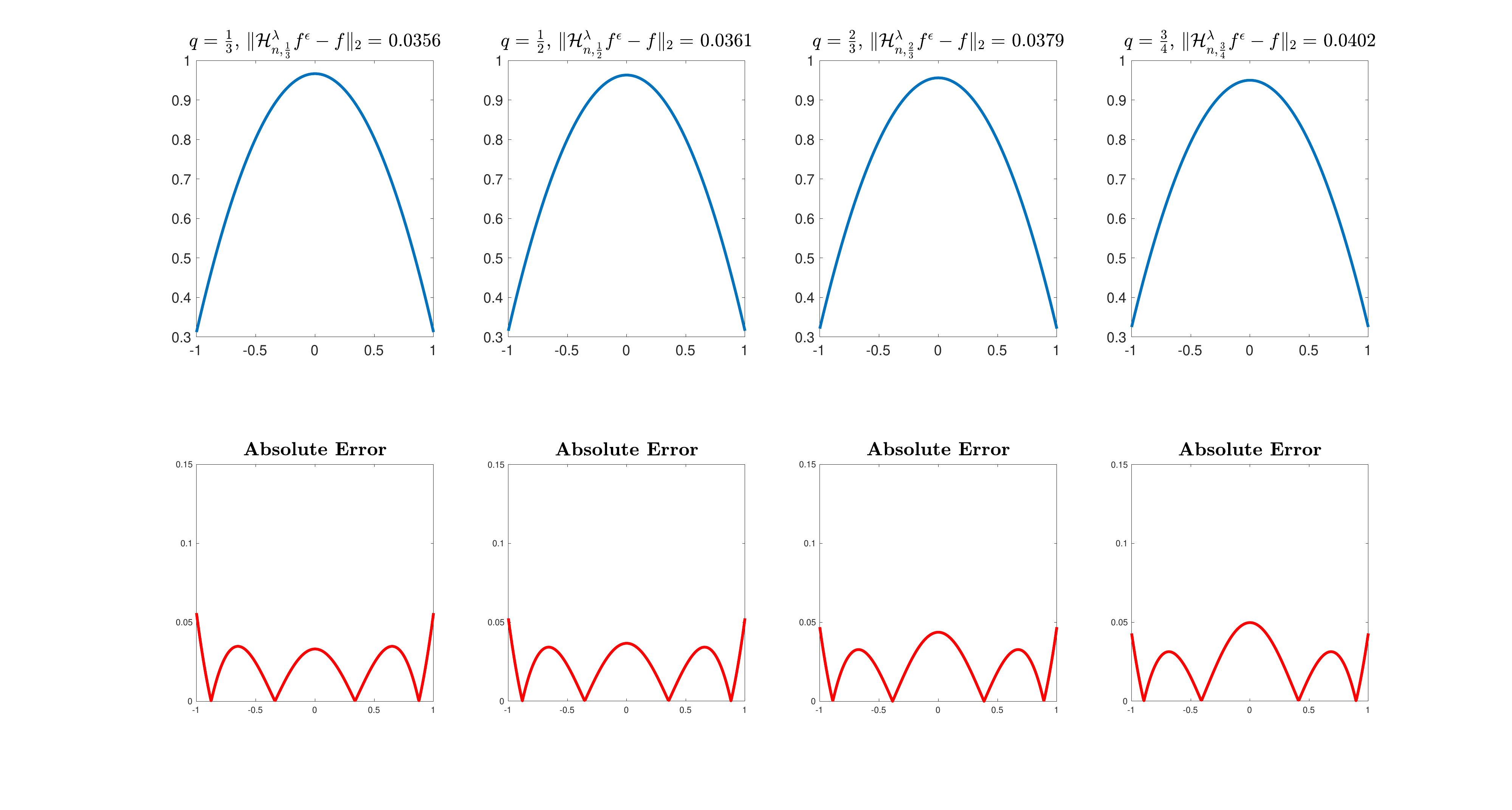}\\
  \caption{Newton hyperinterpolations with $q= \frac{1}{3},\frac{1}{2},\frac{2}{3},\frac{3}{4}$ reconstruct the original function $f=\exp(-x^2)$ from noisy measurements .}\label{fig:interval2}
\end{figure}

Let $f(x)=\exp(-x^2)$ with $x\in[-1,1]$, and $\mathbf{A}\in \mathbb{R}^{N\times d_n}$ with $(N,d_n)=(400,251)$. The sampling sets $\{f(x_j)\}_{j=1}^{N}$ is contaminated by Gaussian noise $\mathcal{N}(0,0.15^2)$ and Impulse noise $\mathcal{I}(0.5)$. The original function $f$ and its reconstructions are shown in Figures \ref{fig:interval} and \ref{fig:interval2}, where the regularization parameter $\lambda$ is chosen to retain only the 2 largest coefficients by magnitude in the hyperinterpolations.

%\subsection{Image processing}

%\subsection{Unit sphere}

\section{Final remark}
This study highlights the significance of recovery thresholding hyperinterpolations in addressing challenges in signal processing, particularly when dealing with noisy data. The proposed nonconvex regularization techniques demonstrate promising capabilities in maintaining sparsity while effectively reconstructing signals. The comparative analysis of various hyperinterpolation methods underscores their applicability and robustness in real-world scenarios, including image processing and function denoising. Future research can further refine these methods, exploring their integration with advanced algorithms to enhance performance in increasingly complex datasets. The findings provide a solid foundation for ongoing exploration in this vital area of mathematical and statistical sciences.

\section{Acknowledgment}
We would like to express our gratitude to Dr. Dingtao Peng and Dr. Hao-Ning Wu for their valuable insights and supportive communications regarding our study.

\newpage
\bibliographystyle{siamplain}
\bibliography{refferences}
\end{document}